\documentclass{article}

\newcommand{\dual}{\mbox{$A^{\ast\ast}$}}

\newcommand{\ro}{\mbox{$\varrho$}}
\newcommand{\C}{\mbox{$C^\ast$}-}

\begin{document}
\begin{center}
{\bf \huge  Orthogonal Pure States in Operator Theory}\\[1cm]
{\large Jan Hamhalter}\\[2cm]
\end{center}
{\small Abstract: We summarize and deepen existing results on
 systems of orthogonal pure states in the context
of JB algebras and \C algebras. Especially, we focus on noncommutative
generalizations of some principles of topology of locally compact spaces
such as exposing points by continuous functions, separating sets by continuous
functions, and multiplicativity of pure states.
}

\section{Introduction and preliminaries}

The aim of this paper is to summarize and deepen recent results on the systems of
orthogonal pure states in operator theory and comment on their relations to axiomatic
foundations of quantum mechanics. There are  two main motivations for our
study of pure states and their orthogonal systems. At first, in the standard
axiomatics of quantum mechanics pure states of a physical system are
identified with the rays in a separable Hilbert space. In this approach
orthogonal pure states, embodying mutually exclusive states of the system,
correspond to orthogonal one-dimensional subspaces. For this reason, the
systems of pure states have been studied in connection with fundamental
questions of quantum mechanics such as hidden variables problem, description
of independence of quantum system, etc. Secondly, from the mathematical
standpoint, pure states encode much of the structure of operator algebras.
Any commutative $C^\ast$-algebra $\cal A$ may be  identified  with
the algebra of complex continuous functions on a locally compact space
$X$, while
pure states  correspond to points of $X$.
There are many important principles on separating points and sets
in topology of locally compact spaces, as for example
the Uryshon Lemma and others. The objective of this note is to overview
possible extensions
of various principles of this kind to the realm of
the \C algebras and the Jordan-Banach algebras.\par
 As the
starting point we recall basic properties of commutative operator algebras.
 Any commutative \C algebra
is isomorphic to the algebra $C_0(X)$ of all continuous  complex functions
on a locally compact Hausdorff space $X$ vanishing at infinity.
 Associative JB algebras are nothing but
algebras of all real functions from $C_0(X)$. At this level there is no
substantial
difference between the Jordan algebras and the \C algebras.
 One of the classical results of Banach function algebras says that
a state \ro{} on $C_0(X)$  is pure if and only if it is multiplicative, i.e.
precisely when
$\ro(fg)=\ro(f)\,\ro(g)$ for all $f,g\in C_0(X)$. Moreover, there is a
one-to-one
correspondence between the pure states and the points of $X$ in the sense that
every pure state \ro{} is of the form $\ro(f)=f(x)$,  ($f\in
C_0(X)$),
 where $x\in X$. These results are far from being true in
the noncommutative case. In fact,  all pure states on a \C algebra are
multiplicative exactly when this algebra  is commutative.
However, we shall establish that even in the noncommutative case
pure states are multiplicative on large
parts of operator algebras.

Another fact concerning exposing point by a continuous function
 will be dealt with in more general
noncommutative setting.  If the compact space $X$ is metrizable,
 then any point $x\in X$ can be
exposed by a continuous function $f$ on $X$ with values in $[0,1]$,
meaning that $x$ is the only point such that $f(x)=1$. Among others we will
show that this principle remains valid for separable algebras.

Finally, working in locally compact space $X$, it is easy to verify that
given a sequence $(x_n)$ discrete as a subspace of $X$,  we can always
find a sequence $(f_n)$ of continuous function on $X$, $0\le f_n\le 1$,
having pairwise disjoint supports and satisfying $f_n(x_m)=\delta_{nm}.$
In other words, points can  be separated by functions. We will investigate to
what extend can this fact be extended to nonassociative generalization of
compact Hausdorff spaces given by JB algebras.

 What makes the problems outlined above
interesting and difficult is the fact that "points" and "sets" cannot be
localized in one space. Namely, in general operator setting points
correspond to elements of  the dual spaces while open sets correspond to
special projection in the double duals.

Let us now recall basic notions and fix the notation.
 Our basic structure will be JB algebra.
A {\em JB algebra} is a real Banach algebra $(A,\circ)$ such that
for all $a,b\in A$,
(i) $a\circ b=b\circ a$,
(ii) $(a\circ b)\circ a^2=a\circ (b\circ a^2)$,
(iii) $||a^2||=||a||^2$,
(iv)  $||a^2+b^2||\geq ||a||^2$.
For all unmentioned details on
 operator algebras we refer to monographs \cite{Davidson, Olsen-Stormer,
 Kadison,Pedersen,Takesaki}.
In what follows $A$ will stand for a JB algebra and $B(H)$ for the algebra
of all bounded operators acting on a Hilbert space $H$.
An important example of a JB algebra is the algebra of all selfadjoint
operators of a \C algebra $\cal A$ endowed with the product $a\circ b= 1/2
(ab+ba)$. A {\em state} on a JB algebra $A$ is a norm one functional $\ro$
on $A$ such that $\ro(a^2)\ge 0$ for all $a\in A$. The {\em pure} state is
an
extreme point of the convex set of all states. Two canonical examples
of pure states  will
be  important for us. At first, the Dirac measures $\delta_x$, $x\in
X$, defined on the function algebra $C_0^R(X)$ of all real continuous
functions on a locally compact  space $X$ vanishing at infinity by
$\delta_x(f)=f(x)$,  $(f\in C_0(X)$,  are the only pure states on
$C_0^R(X)$. Secondly, the vector state on the self-adjoint part of a \C algebra
$\cal A$ acting on a Hilbert space $H$ associated with a unit vector $\xi\in H$
by the formula $\omega_\xi(a)=(a\xi, \xi)$ ($a\in {\cal A}$)  is a pure state
that can be considered as the noncommutative version of the Dirac measure.
 Pure states $\ro_1$ and $\ro_2$
on  $A$ are called {\em orthogonal} if they have the maximal
possible distance,
i.e. if $\|\ro_1-\ro_2\|=2$. The vector states on the algebra of all
self-adjoint operators on a given Hilbert space are orthogonal exactly when the
corresponding vectors are orthogonal. In case of associative algebras all
distinct pure states are orthogonal. It is a key property of the quantum
model that there are distinct but not orthogonal pure states (see Bohr's concept
of complementarity).
\section{Individual pure state}
The problem of exposing point of a compact space by a continuous function
has given impetus to the following definition.
An element $a\in A,$ $||a||\leq 1$, $a\geq 0,$ is {\em determining} for a
pure state $\varrho$ on $A$ if $\varrho(a)=1$ and $\varphi(a)<1$ for any pure
state $\varphi$ on $A$ different from $\varrho$.
In the physical terminology, the state $\ro$ admits the determining element if and
only if there is some observable $O$ with values in [0,1] such that the system
 is in the state \ro{}
precisely when we obtain value 1 when measuring $O$. As an example,
an atomic projection of a Hilbert space $H$ onto one--dimensional subspace
spanned by the unit vector $\xi\in H$ is a determining element for the
vector state $\omega_\xi$ on the algebra of all self-adjoint operators
on a Hilbert space $H$. This is a typical situation in the  standard formalism
of quantum mechanics.
 The following proposition describes states
admitting determining elements:\\

\noindent{\bf 2.1. Proposition.} (\cite{Pure})
 {\em If $A$ is a unital algebra then a pure state $\varrho$ on $A$
admits a determining element if and only if its left ideal
 $L_\varrho =\{ \varrho(x^2)=0\ |\ x\in A\}$
  has a countable approximate unit. Especially, any pure state on a
  separable JB algebra has a determining element.  }\\

The proof is based on the analysis of approximate units and pure states.
This result  is far from being true
for nonseparable algebra. Indeed,  in contrast to Proposition 2.1, only
 normal pure states
on JBW algebras qualify for having determining elements. Let us recall that a JB
algebra is called a JBW algebra if it is a dual Banach space. A state on JBW
algebra is called normal if it can be identified with element from the
predual.\\

\noindent  {\bf 2.2. Theorem.} (\cite{Pure}) {\em A~pure state $\varrho$ on a JBW
  algebra $M$ has a determining  element if and only if $\varrho$ is
  normal.
  Especially, if $M$ is a von Neumann algebra then $\varrho$ has a
  determining element if and only if it is a vector state concentrated on
  a direct summand isomorphic to $B(H)$. }\\

The previous Theorem says that the states described by the rays in a Hilbert space $H$ are the
only pure states on von Neumann  algebras having determining elements.
 This  advocates  basic assumptions in early models of quantum
mechanics.

\section{Orthogonal pure states and supporting systems}

In this section we shall deal with a more complicated case of a system of
orthogonal pure states. Suppose that $(\varrho_n)$ is a sequence of
(pairwise)  orthogonal pure states on $A$.
  A~sequence $(b_n)$ in $A$ with $0\leq b_n\leq 1$, and $b_n\circ b_m=0$
  whenever $n\not=m$ is called {\em supporting} for $(\varrho_n)$ if
\[\varrho_n(b_m)=\delta_{mn}\qquad \mbox{ for all } n,m\,. \]
  If, in addition, each $b_n$ is determining for $\varrho_n,$ we call
  $(b_n)$ a {\em determining supporting sequence}.
In example of associative algebra in the introduction we have seen that
a sequence of orthogonal pure states has a supporting sequence if the
corresponding points can be separated by open sets. In order to find an
analogy of this result for nonassociative algebras we shall need the
following concepts. Let $A$ be canonically embedded into its second dual
$A^{\ast\ast}.$ It is known that \dual{} is a JBW algebra extending the product in $A$.
 A~projection $p$ in $A^{\ast\ast}$ is called {\em open} if there is an
 increasing net of elements $(a_\lambda)$ in $A$ with $a_\lambda\nearrow
 p$.
A~projection $p$ in $A^{\ast\ast}$ is {\em closed} if $1-p$ is open.
Any  state $\ro$ on $A$ canonically extends to a normal state on
$A^{\ast\ast}$. Hence, there is a smallest projection $s(\ro)$ in \dual{}
such that  $\ro(s(\ro))=1$. We call $s(\ro)$ the {\em
support projection} of $\ro$. It is known that $s(\ro)$ is always
minimal and closed if \ro{} is pure. It was proved in \cite{Supporting} that a system $\ro_\alpha$
of orthogonal pure states admits a supporting system if and only if there
is a system $(p_\alpha)$ of open, pairwise orthogonal projections in
\dual{}, such that $s(\ro_\alpha)\le p_\alpha$ for each $\alpha$.
 This is a precise analogy of the classical case.    Our main result gives
 sufficient conditions for the existence of supporting system in terms of
 primitive ideal space. Let $P(A)$ be the set of all pure states on $A$
 equipped with the weak$^\ast$-topology. Let $c(\ro)$ be the {\em central
 cover} of $\ro$, i.e. the smallest central projection in \dual{} majorizing the
 support projection $s(\ro)$. We  define a homomorphism $\pi_{\ro}: A\to
 \dual$ by $\pi_{\ro}(x)=c(\ro)\circ x$, $x\in A$. In case of \C algebras
 $\pi_{\ro}$ corresponds to the G.N.S. representation of $\ro$. The {\em
 primitive ideal space} is now defined as  the set of all kernels
 of the representations, $Prim(A)= \{ Ker\, \pi_{\ro}\ |\ \ro\in P(A)\}$. The
 primitive ideal space is endowed with the Jacobson topology whose closure
 operation is $S\subset Prim(A)\to \overline{S}=\{F\in Prim(A)\ |\ F\supset
 \cap S\}$. The canonical map $\tau$ now sends pure state \ro{} to its
 kernel $Ker \ro$. A set $P$ in the pure state space $P(A)$ is called {\em almost
 separated} if its image $\tau(P)$ in $Prim(A)$ can be covered by disjoint open sets
 $(U_\alpha)$  such that each open set $\tau^{-1}(U_\alpha)$
 contains finitely many elements. The following result says that any
 almost separated sequence of orthogonal pure states admits a supporting
 sequence.\\

\noindent  {\bf 3.1. Theorem.} (\cite{Supporting}) {\em\   Any sequence $(\varrho_n)$ of
   almost separated orthogonal
 pure states on a JB algebra $A$ has a supporting sequence. This sequence can be chosen
 determining if $A$ is separable. }\\

 By combining the Hahn-Banach and the Krein-Milman theorems we can derive the
 following consequences of the previous results for the restriction properties
 of states. By Theorem 3.1 any finite sequence $\ro_1,\ro_2,\ldots \ro_n$ of
 orthogonal pure states on a separable algebra $A$ has a determining sequence $b_1, b_2,\ldots,
 b_n$. The algebra $C$ generated by $b_1,b_2,\ldots, b_n$ is associative and all
 states $\ro_1,\ro_2,\ldots,\ro_n$ restrict simultaneously to pure states
 on $C$ and are uniquely given by these restrictions. If $n=2$ then the
 algebra $C$ can even be specified to be singly generated. From the physical point
 of view, two mutually exclusive states of the system are given by the structure of
 just one observable. Since any associative subalgebra can  be extended to
 maximal one, we infer that any sequence $(\ro_n)$ of almost separated states on
 a separable algebra admit a {\em determining} maximal associative subalgebra
 in the sense that each pure state $\ro_n$ is uniquely determined by its
 restriction to this algebra.

  The situation concerning infinitely many orthogonal pure states is
  discussed in the following Theorem. Following \cite{A.A.P.1} we say that a
     sequence $(\ro_n)$ of pure states on a JB algebra $A$ {\em approaches}
     to infinity if $\lim_{n\to\infty}\ro_n(a)=0$ for every $a\in A$ such
     that the spectrum of $a$ contains zero. For example, if $A$ is the real
     part of the algebra of all continuous functions on the real line
     vanishing at infinity, then $(\delta_{x_n})$ approaches to infinity if
     and only if $|x_n|\to \infty$.\\

\noindent  {\bf 3.2. Theorem.} {\em Let $A$ be a separable JB algebra and $(\ro_n)$
  be a sequence of orthogonal pure states approaching to infinity. The
  following statements are equivalent:
  \begin{enumerate}
  \item[{\rm (i)}] $(\ro_n)$ has a supporting system.
  \item[{\rm (ii)}]  $(\ro_n)$ has a determining supporting system.
  \item[{\rm (iii)}]  There exists a maximal associative subalgebra $C$ of
  $A$ determining for $(\ro_n)$.
  \item[{\rm (iv)}]  There exists a maximal associative subalgebra $C$ of
  $A$ such that  $(\ro_n|C)$ forms an orthogonal sequence of pure states.
    \end{enumerate}}

    Proof: The equivalence of (i) and (ii) is the content of
    \cite[Proposition 3.1]{Supporting}\par
    (ii)$\Rightarrow$(iii) Suppose  that $(b_n)$ is a supporting determining
    sequence for $(\ro_n)$. The algebra $C'$ generated by $b_n's$ is
    associative. Let $C$ be any maximal associative subalgebra extending
    $C$. Any state $\ro_n|C'$ is pure and extends uniquely to a pure state
    on $A$. Hence, state $\ro_n|C'$ has a unique extension to a state
    $\ro_n|C$ which has to be, by the Krein-Milman theorem, a pure state.
    Hence, each $\ro_n$ is uniquely  determined by its pure restriction to
    $C$. \par
    (iii)$\Rightarrow$ (iv) Let $C$ be an associative subalgebra
    fulfilling (iii). Since $\ro_n$'s have to be distinct on $C$ we have
    $\| (\ro_n-\ro_m)| C\|=2$ because $C$ is associative.\par
    (iv)$\Rightarrow$ (i)
     Applying Theorem 3.1 to a separable subalgebra $C$ we get a desired
     supporting sequence for $(\ro_n)$. The proof is completed.\\

  These extension results contribute to the line of research in works
   \cite{Akemann1,A.A.P.2,Anderson,Barnes,Bunce}.

  \section{Multipicativity of orthogonal pure states}

  As we have seen the pure state on an operator algebra is seldom
  multiplicative. On the other hand, it can be uniquely
  determined by multiplicative state on some maximal associative subalgebra
  provided that the algebra is separable. As it was shown by Akemann
  this does not hold for nonseparable algebras \cite{Akemann1}. However, we can
  still ask whether there exists some maximal associative subalgebra such that
  given pure state is multiplicative on it. The following problem, that is
  central for this section, seems to be a natural operator-theoretic
  extension of the
  equivalence between pureness and multiplicativity for states on function
  algebras.\\

\noindent  {\bf 4.1. Problem.} {\em Is a given sequence of orthogonal pure states on a JB
  algebra (resp. \C algebra) pure (i.e. multiplicative) on some maximal
  associative (resp. abelian) subalgebra? }\\

  There are many partial positive results along this line.
  E.Stormer haracterized pure states in terms of definite subasets
  \cite{Stor}.
  Aarnes and
Kadison answered the problem in the positive for a pure state on a
  separable unital \C algebra \cite{A.K.}. Then Akemann extended their result
  to finitely many states on a separable (not necessarily) unital \C
  algebra \cite{Akemann1}. It was Barnes who realized that what is important is the
  separability of the G.N.S. representation of a given state.
   He established the positive answer for
   individual pure state with separable representation \cite{Barnes}.
   Also he gave, independently on Bunce \cite{Bunce},
   positive  answer for pure states on Type I \C algebras \cite{Barnes}. The topic was then revived
  by Akemann, Anderson and Pedersen in the paper \cite{A.A.P.2} dealing with
  states approaching to infinity. They have proved that pure orthogonal
  sequence of nearly inequivalent states on \C algebra
   approaching to infinity  restrict simultaneously to multiplicative states on some
  maximal abelian subalgebra. Unfortunately, they needed additional
  assumption on accumulation points of this sequence. Despite these positive
  results we have proved that the answer to Problem 4.1 is in the negative.
  Moreover, it turns out that the counterexample is quite generic.\\

\noindent  {\bf 4.2. Counterexample.} (\cite{Multiplicativity}) {\em Let $\cal A$ be a separable unital infinitely
  dimensional \C algebra acting irreducibly on a Hilbert space $H$.
  Let $K$ be the algebra of all compact operators acting on $H$. Suppose
  that ${\cal A}/ {\cal A}\cap K$ is noncommutative. Then there exists a sequence of
  orthogonal pure states on $\cal A$ that do not restrict (simultaneously)
  to pure states on any maximal commutative \C subalgebra of $\cal A$. }\\

  Proof: The irreducible algebra $\cal A$ either contains all compact operators
  or ${\cal A}\cap K=\{0\}$. The proof of  both cases is given
  in \cite{Multiplicativity}. Here we present a simpler version of arguments
   concerning the
  case when ${\cal A}\cap K=\{0\}$. Since $H$ is separable we can find a dense
  sequence $(x_n)$ in its unit sphere.  It means that the convex hull
  of the sequence $(\omega_{x_n})$ of the corresponding vector states is weak$^\ast$-dense in
  the state space of $A$ (see e.g. \cite{Kadison}). Set
  \[ f=\sum_n \frac 1{2^n}\omega_{x_n}\,. \]
   By the improved Glimm's lemma
  \cite{A.A.P.2,Glimm} there is an orthonormal basis $(\xi_n)$ of $H$ such that
  \[ f(a)=\lim_{n\to\infty}(a\xi_n, \xi_n)\,. \]
   By taking the corresponding pure states $(\omega_{\xi_n})$ we obtain a sequence
   of pure orthogonal states on $\cal A$ ($\cal A$ is irreducible). Suppose that $C$ is a maximal abelian
   subalgebra such that all states $(\omega_{\xi_n})$ are multiplicative on
   $C$. Then, obviously,
    $f$ is also multiplicative on $C$. In other words
   $f$ is pure on $C$ and so $f=\omega_{\xi_n}$ for all $n$ on $C$. By the density
   of the convex hull of $\omega_{x_n}$, we see that all states on $C$ are
   multiples of $f$. By the spectral theory $C$ has to be one-dimensional
   and so $\dim {\cal A}=1$ -- a contradiction with irreducibility of ${\cal A}$ on an
   infinite-dimensional $H$.\par
      The case when ${\cal A}$ contains the ideal of compact operators is based on
    the Weyl-von Neumann theorem and may be found in \cite{Multiplicativity}.\\

    This counterexample has interesting consequence for the Calcin algebra,
    i.e. for the quotient of $B(H)$ by the ideal of compact operators $K$. Namely,
    we can construct a separable subalgebra of the Calcin algebra having
    a sequence of orthogonal pure states that are not simultaneously
    multiplicative on any maximal commutative $C^\ast$-algebra. In contrast
    to this, Anderson proved in \cite{Anderson-Calcin} that for any sequence of states
    on the Calcin algebra there is a maximal associative subalgebra such
    that all members of this sequence are multiplicative on it. This
    indicates that the interplay between the size of \C-algebra and the size
    of its maximal abelian subalgebras is quite delicate.\par
     Counterexample 4.2 indicates possible limitations to improvements of
     positive results. We have to put some additional conditions on
     the system of orthogonal pure states as it is shown in the following
     theorem.  \\

\noindent     {\bf 4.3. Theorem.} (\cite{Multiplicativity}) {\em Let $(\varrho_n)$ be a sequence of orthogonal pure
states on a JB  algebra $A$ approaching to infinity and such that
$c(\varrho_n)A^{\ast\ast}$ is $\sigma$-finite for all $n$.
 Suppose further that $\sum_{n=m}^\infty s(\varrho_n)$ is a closed
 projection for all $m$.
 Then there is a
maximal associative subalgebra $B$ of $A$ such that all states
$(\varrho_n)$ restrict to pure states on $B$. \par
In particular, for a finite system of orthogonal pure states on $A$ with
$\sigma$-finite central covers there is always a maximal associative
 subalgebra such that all states are pure on it}.\\

 This result is sharper than existing results even in case of \C algebras \cite{A.K.,
  Akemann1, A.A.P.1, A.A.P.2,Anderson,Barnes,Bunce}.
  However,
 the main difficulty lies in the Jordan case when we have to develop new
 methods to overcame the lack of Hilbert-space representations which are
 ubiquitous in \C case.\par
  Summing it up, the question on multiplicativity of an infinite sequence of
  orthogonal pure states is more or less solved. Nevertheless, the following well
  known question is still open: Does every pure state on a JB algebra
  restrict to a pure state on some maximal associative subalgebra?\\

  The results on the determining elements and multiplicativity of
  orthogonal pure states have been applied to studying the structure of
  the compact JB algebras \cite{Supporting} and can  also be of some interest for the pure
  extension property studied e.g. in \cite{Bunce1}.

\vspace{5cm}

\noindent{\bf Acknowledgement}: The  author would like to express his gratitude
to the Alexander von Humboldt Foundation, Bonn, for supporting his
research. He also acknowledges the support of the Grant Agency of the Czech
Republic, Grant no. 201/03/0455, "Noncommutative measure theory",
 and the support of the
 Czech Technical University, Grant. No. MSM 210000010, "Applied
 Mathematics in Technical Sciences".

Jan Hamhalter, Czech Technical University -- El. Eng, Department of
Mathematics, Technicka 2, 166 27, Prague 6, Czech Republic.
e-mail: hamhalte@math.feld.cvut.cz

\end{document}